# Monotonicity, asymptotic normality and vertex degrees in random graphs

SVANTE JANSON

[1]*Department of Mathematics, Uppsala University, PO Box 480, SE-751 06 Uppsala, Sweden.
E-mail: svante.janson@math.uu.se. url: http://www.math.uu.se/~svante/*

We exploit a result by Nerman which shows that conditional limit theorems hold when a certain monotonicity condition is satisfied. Our main result is an application to vertex degrees in random graphs, where we obtain asymptotic normality for the number of vertices with a given degree in the random graph $G(n,m)$ with a fixed number of edges from the corresponding result for the random graph $G(n,p)$ with independent edges. We also give some simple applications to random allocations and to spacings. Finally, inspired by these results, but logically independent of them, we investigate whether a one-sided version of the Cramér–Wold theorem holds. We show that such a version holds under a weak supplementary condition, but not without it.

## 1. Introduction

Many random variables in different areas of probability, statistics and combinatorics can be expressed as some "simpler" random variable conditioned on a specific value of another. A few examples are given in Sections 3 and 4 below; many other can be found in the references.

Such representations are, among other things, useful for the derivation of asymptotic results. Generally speaking, if $(X_n, Y_n) \xrightarrow{d} (X, Y)$, we would like to conclude that the conditional distributions also converge, that is, $(X_n \mid Y_n = y) \xrightarrow{d} (X \mid Y = y)$. Of course, this is not true in general, but it holds in many cases and several authors have proved more or less general theorems of this type under various assumptions; see, for example, Steck [28], Holst [9, 11], Janson [17] and the further references given there. (For a trivial counterexample which shows that some assumptions are needed, let $\mathbb{P}(X_n = Y_n = 0) = 1/n$ and $\mathbb{P}(X_n = 1, Y_n = 1/n) = 1 - 1/n$, with $X = 1$ and $Y = 0$. Then, $(X_n \mid Y_n = 0) = 0$, but $(X \mid Y = 0) = 1$ a.s.)

The purpose of the present paper is to exploit a result of Nerman [23] which shows that such conditional limit theorems hold in the special, but not uncommon, situation that a certain monotonicity condition holds. This result seems to have been somewhat neglected, but it has many applications. We illustrate its power first by some simple applications to random allocations. Our main result is an application to vertex degrees in random graphs, where we derive a new result for the random graph $G(n, m)$ with a fixed number of edges from the corresponding (but slightly different) result for the random graph $G(n, p)$ with independent edges; indeed, it was this problem that led to the present paper. The method also applies to other properties of $G(n, m)$ and we mention some of these.

We state Nerman's theorem, in versions suitable for easy applications, in Section 2. In Section 3, we illustrate the theorems by a simple application to random allocations, where we give short proofs of some known results. We present our main application to random graphs in Section 4, where we state and prove Theorem 4.1 on vertex degrees. This also illustrates how results for other monotone functions of random graphs may be obtained. An application of Theorem 4.1 to the study of the $k$-core is given in [18]; this application was the original motivation for the research that led to the present paper.

We provide another simple application, to spacings, in Section 5.

It often happens that vector-valued versions of limit theorems follow from the one-dimensional versions by the Cramér–Wold device, that is, by considering linear combinations of the components; see Cramér and Wold [6] or, for example, Billingsley [4], Theorem 7.7. This is not the case here, since the assumptions in the multidimensional case of Theorem 2.2 allow us to apply the one-dimensional case only to linear combinations $t_1 X_n^{(1)} + \cdots + t_d X_n^{(d)}$, where all $t_i$ have the same sign. Although not needed for our results (the multidimensional case is proved by Nerman as the one-dimensional case), we find it interesting to investigate whether, in general, it is enough to show convergence for such linear combinations, that is, whether there is a one-sided version of the Cramér–Wold device. In Section 6, we show that the answer is affirmative under a weak supplementary condition, but not in general. This is equivalent to a corresponding uniqueness problem and to the question of whether a characteristic function is determined by its restriction to the first octant.

All unspecified limits below are as $n \to \infty$. Further, as usual, $\delta_{ij}$ is 1 when $i = j$ and 0 otherwise.

## 2. Nerman's general results

***Definition.*** *Let $X$ and $Y$ be random variables defined on the same probability space. We say that $X$ is* stochastically increasing *with respect to $Y$ if the conditional distributions $\mathcal{L}(X \mid Y = y)$ are stochastically increasing in $y$, that is, if*

$$\mathbb{P}(X \leq x \mid Y = y_1) \geq \mathbb{P}(X \leq x \mid Y = y_2) \qquad \textit{for any real } x \textit{ and } y_1 \leq y_2. \tag{2.1}$$

*If $Y$ has a discrete distribution, we may, here and below, consider only $y$ (and $y_1, y_2$) such that $\mathbb{P}(Y = y) > 0$ and there is no problem with defining the conditional distributions* *and probabilities. In general, for example, for continuous $Y$, the conditional distribution $\mathcal{L}(X \mid Y = y)$ is defined only up to equivalence, that is, for a.e. $y$ with respect to the distribution $\mathcal{L}(Y)$. The precise definition is that there exists a version of $y \mapsto \mathcal{L}(X \mid Y = y)$ that is stochastically increasing in $y$; it is this version that is used below. (Thus, (2.1) holds with the conditional probabilities defined by this version.)*

*We say that $X$ is* stochastically decreasing *with respect to $Y$ if $-X$ is stochastically increasing, and* stochastically monotone *with respect to $Y$ if it is either stochastically increasing or stochastically decreasing with respect to $Y$.*

The definitions extend to vector-valued $X$ and $Y$, using the partial order on $\mathbb{R}^d$ defined by $(s_1, \ldots, s_d) \leq (t_1, \ldots, t_d)$ if $s_i \leq t_i$ for each $i$.

***Remark 2.1.*** *It is well known that if $X$ is real-valued, then (2.1) is equivalent to the existence of an increasing coupling of $(X \mid Y = y_1)$ and $(X \mid Y = y_2)$, that is, a pair of random variables $\widetilde{X}_1$ and $\widetilde{X}_2$ such that $\widetilde{X}_j \stackrel{d}{=} (X \mid Y = y_j)$, $j = 0, 1$, and $\widetilde{X}_1 \leq \widetilde{X}_2$ a.s. This is not generally true in the vector-valued case, but the (perhaps more natural) condition that there always exists such an increasing coupling is stronger and implies (2.1).*

We can now state Nerman's theorem. Let $\mathcal{P}(\mathbb{R}^q)$ denote the set of probability measures on $\mathbb{R}^q$, equipped with the usual (weak) topology.

***Theorem 2.2.*** *Suppose that $(X_n, Y_n)$, $n \geq 1$, are pairs of random vectors, with $X_n \in \mathbb{R}^q$ and $Y_n \in \mathbb{R}^r$ for some $q, r \geq 1$, such that $X_n$ is stochastically monotone with respect to $Y_n$. Further, suppose that for some sequences of real numbers and vectors $a_n > 0$, $b_n \in \mathbb{R}^q$, $c_n > 0$, $d_n \in \mathbb{R}^r$,*

$$(a_n^{-1}(X_n - b_n), c_n^{-1}(Y_n - d_n)) \xrightarrow{d} (X, Y)$$

*for a pair of random vectors $X \in \mathbb{R}^q$ and $Y \in \mathbb{R}^r$. Assume that $y_n$ is a sequence in $\mathbb{R}^r$ such that $c_n^{-1}(y_n - d_n) \to \xi$ for some $\xi \in \mathbb{R}^r$ and let $\widetilde{X}_n$ be a random vector whose distribution equals the conditioned distribution $\mathcal{L}(X_n \mid Y_n = y_n)$. Finally, suppose that $\xi$ is an interior point of the support of $Y$ and that there exists a version of $y \mapsto \mathcal{L}(X \mid Y = y)$ that is continuous at $y = \xi$ as a function of $y \in \mathbb{R}^r$ into $\mathcal{P}(\mathbb{R}^q)$. Then,*

$$a_n^{-1}(\widetilde{X}_n - b_n) \xrightarrow{d} \mathcal{L}(X \mid Y = \xi).$$

**Proof.** Nerman [23], Theorem 1 and Remark, proved the case $a_n = c_n = 1$, $b_n = 0$, $d_n = 0$. The general version follows immediately by replacing $(X_n, Y_n)$ with $(a_n^{-1}(X_n - b_n), c_n^{-1}(Y_n - d_n))$. □

The case when $X$ and $Y$ have a joint normal distribution is perhaps the most interesting, both because it appears in many applications and because the result simplifies somewhat. In this case, assuming that the covariance matrix of $Y$ is non-singular, it is elementary and well known that there exists a continuous version of $y \mapsto \mathcal{L}(X \mid Y = y)$, given by

$$(X \mid Y = y) \stackrel{d}{=} X + A(y - Y), \tag{2.2}$$

where $A$ is the $q \times r$ matrix given by $A := \operatorname{Cov}(X,Y)(\operatorname{Var}(Y))^{-1}$. Note that this, too, has a normal distribution. (To see (2.2), note that $Z := X - AY$ and $Y$ are uncorrelated and thus independent. Since $X = Z + AY$, it follows that $(X \mid Y = y) \stackrel{d}{=} Z + Ay = X + A(y - Y)$.)

For ease of application, we state the result in the normal case separately, restricting ourselves to the case $r = 1$, which simplifies notation and is the most important case for applications.

**Theorem 2.3.** *Suppose that $(X_n, Y_n)$, $n \geq 1$, are pairs of random vectors $X_n = (X_n^{(1)}, \ldots, X_n^{(q)}) \in \mathbb{R}^q$ and variables $Y_n \in \mathbb{R}$ for some $q \geq 1$, such that $X_n$ is stochastically monotone with respect to $Y_n$. Suppose, further, that for some sequences of real numbers and vectors $a_n > 0$, $b_n \in \mathbb{R}^q$, $c_n > 0$, $d_n \in \mathbb{R}$,*

$$(a_n^{-1}(X_n - b_n), c_n^{-1}(Y_n - d_n)) \stackrel{d}{\longrightarrow} (X, Y) \tag{2.3}$$

*for a normally distributed random vector $(X, Y)$ with $X = (X^{(1)}, \ldots, X^{(q)}) \in \mathbb{R}^q$ such that $\operatorname{Var}(Y) > 0$. (Thus, $X^{(1)}, \ldots, X^{(q)}, Y$ are jointly normal.)*

*Assume that $y_n$ is a sequence in $\mathbb{R}$ such that $c_n^{-1}(y_n - d_n) \to \xi$ for some real $\xi$ and let $\widetilde{X}_n = (\widetilde{X}_n^{(1)}, \ldots, \widetilde{X}_n^{(q)})$ be a random vector whose distribution equals the conditioned distribution $\mathcal{L}(X_n \mid Y_n = y_n)$.*

*Then, with $\gamma = (\gamma^{(1)}, \ldots, \gamma^{(q)})$, where $\gamma^{(i)} = \operatorname{Cov}(X^{(i)}, Y)/\operatorname{Var}(Y)$,*

$$a_n^{-1}(\widetilde{X}_n - b_n) \stackrel{d}{\longrightarrow} \widetilde{X} := X + (\xi - Y)\gamma.$$

*Thus, $\widetilde{X} = (\widetilde{X}^{(1)}, \ldots, \widetilde{X}^{(q)})$ is normal with*

$$\mathbb{E}\widetilde{X} = \mathbb{E}X + (\xi - \mathbb{E}Y)\gamma,$$
$$\operatorname{Cov}(\widetilde{X}^{(i)}, \widetilde{X}^{(j)}) = \operatorname{Cov}(X^{(i)}, X^{(j)}) - \operatorname{Cov}(X^{(i)}, Y)\operatorname{Cov}(X^{(j)}, Y)/\operatorname{Var}(Y).$$

In the one-dimensional case, this result may be stated as follows.

**Corollary 2.4.** *Suppose that the assumptions of Theorem 2.3 hold with $q = 1$ and thus $X_n$ and $X$ real-valued. Let $\sigma_X^2 := \operatorname{Var}(X)$, $\sigma_Y^2 := \operatorname{Var}(Y)$, $\sigma_{XY} := \operatorname{Cov}(X,Y)$ and $\gamma := \sigma_{XY}/\sigma_Y^2$. Then,*

$$a_n^{-1}(\widetilde{X}_n - b_n) \stackrel{d}{\longrightarrow} N(\mathbb{E}X + \gamma(\xi - \mathbb{E}Y), \sigma_X^2 - \sigma_{XY}^2/\sigma_Y^2). \tag{2.4}$$

*If, further, $\sigma_X^2 > 0$, then the asymptotic variance in (2.4) equals $(1 - \rho^2)\sigma_X^2$, where $\rho := \sigma_{XY}/(\sigma_X \sigma_Y)$ is the correlation between $X$ and $Y$.*

The variance in (2.4) is the same as the residual variance in linear regression. This extends to the multi-dimensional case.

We can weaken the assumptions in the multi-dimensional case somewhat.

**Corollary 2.5.** *In Theorems 2.2 and 2.3, we can replace the assumption that $X_n$ is stochastically monotone with respect to $Y_n$ by the assumption that $TX_n$ is stochastically monotone with respect to $Y_n$ for some invertible linear operator $T$ on $\mathbb{R}^d$.*

**Proof.** We can apply the theorems to $(TX_n, Y_n)$, with $X$ and $b_n$ replaced by $TX$ and $Tb_n$, respectively. The result follows by applying $T^{-1}$. □

We end this section by stating a companion result of Nerman [23], Section 4, on moment convergence.

**Theorem 2.6.** *If the pth absolute moments of the components of $a_n^{-1}(X_n - b_n)$ converge to the corresponding moments of $X$ in one of the theorems or corollaries above, then all (mixed) moments and absolute moments of order at most $p$ of $a_n^{-1}(\widetilde{X}_n - b_n)$ converge to the corresponding moments of $\widetilde{X}$. In particular, if the means and (co)variances $a_n^{-1}(X_n - b_n)$ converge to those of $X$, then the means and (co)variances of $a_n^{-1}(\widetilde{X}_n - b_n)$ converge to those of $\widetilde{X}$.*

## 3. A simple application: random allocation

*Example 3.1.* Let $m$ balls be thrown independently of each other into $n$ boxes, with probability $1/n$ of landing in each box. Let $N_k$ be the number of balls landing in box $k$, $k = 1, \ldots, n$. (Thus, $(N_1, \ldots, N_n)$ has a multinomial distribution.) Let $Z_{mn}$ be the number of empty boxes. We are interested in asymptotics as $n \to \infty$ and $m = m(n) \to \infty$. This (and various extensions) has been studied by many authors; see, for example, von Mises [22], Feller [7], Section IV.2, Arfwedson [1], Weiss [29], Rényi [25], Rosén [26, 27], Holst [12], Hwang and Janson [14] and the monograph by Kolchin, Sevast'yanov and Chistyakov [20].

In order to apply the results above, we instead throw a random number $M \sim \text{Po}(\lambda_n n)$ balls. The numbers $N_k$ of balls in the different boxes are then i.i.d. with $N_k \sim \text{Po}(\lambda_n)$. (See Holst [12, 13] for similar uses of Poissonization in this and related problems.) We let $X_n$ denote the number of empty boxes, let $Y_n := M$ and observe that

$$(X_n, Y_n) = \sum_{k=1}^{n} (\mathbf{1}[N_k = 0], N_k).$$

The terms in the sum are i.i.d. random vectors with mean $(e^{-\lambda_n}, \lambda_n)$ and covariance matrix given by

$$\text{Var}(\mathbf{1}[N_k = 0]) = e^{-\lambda_n}(1 - e^{-\lambda_n}),$$

$$\mathrm{Var}(N_k) = \lambda_n,$$
$$\mathrm{Cov}(\mathbf{1}[N_k = 0], N_k) = -\mathrm{e}^{-\lambda_n}\lambda_n.$$

It follows from the central limit theorem that if $n \to \infty$ and $\lambda_n \to \lambda > 0$, then

$$n^{-1/2}(X_n - n\mathrm{e}^{-\lambda_n}, Y_n - n\lambda_n) \xrightarrow{d} (X, Y), \qquad (3.1)$$

with $X$ and $Y$ jointly normal with mean 0 and (co)variances $\sigma_X^2 = \mathrm{e}^{-\lambda}(1 - \mathrm{e}^{-\lambda})$, $\sigma_Y^2 = \lambda$, $\sigma_{XY} = -\lambda\mathrm{e}^{-\lambda}$.

It is obvious that $X_n$ is stochastically decreasing with respect to $Y_n = M$, since throwing another ball can only decrease the number of empty boxes. Moreover, if we condition on $Y_n = m$, we are back in the situation of throwing a given number $m$ balls and thus $\tilde{X}_n = Z_{mn}$. As is well known in this and many related situations, we can here take any $\lambda_n > 0$, but for the continuation of the argument, the choice matters. We use the (natural) choice $\lambda_n = m(n)/n$.

Consequently, if $m = m(n)$ is such that $m(n)/n \to \lambda > 0$, then by Corollary 2.4, using $a_n = c_n = n^{1/2}$, $b_n = n\mathrm{e}^{-\lambda_n}$, $d_n = n\lambda_n = m(n)$, $y_n = m(n)$, $\xi = 0$ and $\gamma = \sigma_{XY}/\sigma_Y^2 = -\mathrm{e}^{-\lambda}$, we have

$$n^{-1/2}(Z_{mn} - n\mathrm{e}^{-m/n})$$
$$\xrightarrow{d} X - \gamma Y = X + \mathrm{e}^{-\lambda}Y \sim N(0, \mathrm{e}^{-\lambda} - \mathrm{e}^{-2\lambda} - \lambda\mathrm{e}^{-2\lambda}),$$

as shown by Weiss [29]; see also Rényi [25] and Kolchin *et al.* [20], Theorem I.3.1.

Corollary 2.4 does not apply directly to the number of boxes with exactly one ball, since this is not stochastically monotone with respect to the number of balls. However, denoting this number by $Z_{mn}^{(1)}$, the sum $Z_{mn} + Z_{mn}^{(1)}$ is the number of boxes with at most one ball, which is stochastically decreasing with respect to $m$. Consequently, we can argue, as above, for the vector $(Z_{mn}, Z_{mn} + Z_{mn}^{(1)})$, using Theorem 2.3, and conclude joint asymptotic normality for $(Z_{mn}, Z_{mn} + Z_{mn}^{(1)})$ and thus for $(Z_{mn}, Z_{mn}^{(1)})$. This is a simple instance of Corollary 2.5.

More generally, if $Z_{mn}^{(j)}$ is the number of boxes with exactly $j$ balls, then Corollary 2.5 applies to $(Z_{mn}, Z_{mn}^{(1)}, \ldots, Z_{mn}^{(J)})$, for any fixed $J$, with

$$T(z_0, z_1, \ldots, z_J) = (z_0, z_0 + z_1, \ldots, z_0 + \cdots + z_J). \qquad (3.2)$$

We again assume that $m = m(n)$ is such that $\lambda_n := m/n \to \lambda > 0$ and denote the Poisson probabilities by

$$\pi_\lambda(k) := \mathbb{P}(\mathrm{Po}(\lambda) = k) = \frac{\lambda^k}{k!}\mathrm{e}^{-\lambda}. \qquad (3.3)$$

We can take $a_n = c_n = n^{1/2}$, $b_n^{(j)} = \pi_{\lambda_n}(j)n$, $d_n = y_n = m(n)$ and then (2.3) holds by the central limit theorem; in this case, we have, if we let $W$ denote a random variable with

$W \sim \mathrm{Po}(\lambda)$,

$$\mathrm{Cov}(X^{(i)}, X^{(j)}) = \mathrm{Cov}(\mathbf{1}[W = i], \mathbf{1}[W = j]) = \delta_{ij}\pi_\lambda(i) - \pi_\lambda(i)\pi_\lambda(j),$$
$$\mathrm{Cov}(X^{(i)}, Y) = \mathrm{Cov}(\mathbf{1}[W = i], W) = i\pi_\lambda(i) - \lambda\pi_\lambda(i),$$
$$\mathrm{Var}\, Y = \mathrm{Var}\, W = \lambda$$

and thus

$$\mathrm{Cov}(\widetilde{X}^{(i)}, \widetilde{X}^{(j)}) = \sigma_{ij}^* := \mathrm{Cov}(\mathbf{1}[W = i], \mathbf{1}[W = j])$$
$$- \mathrm{Cov}(\mathbf{1}[W = i], W)\mathrm{Cov}(\mathbf{1}[W = j], W)/\mathrm{Var}\, W \quad (3.4)$$
$$= \delta_{ij}\pi_\lambda(i) - \pi_\lambda(i)\pi_\lambda(j)\left(1 + \frac{(i - \lambda)(j - \lambda)}{\lambda}\right).$$

In other words, jointly for all $j \geq 0$,

$$n^{-1/2}(Z_{mn}^{(j)} - n\pi_{\lambda_n}(j)) \stackrel{d}{\longrightarrow} \widetilde{X}^{(j)},$$

where $\widetilde{X}^{(j)}$, $j \geq 0$, are jointly Gaussian with means 0 and covariances given by (3.4), as shown by other methods in Kolchin *et al.* [20], Theorem II.2.3; see also Békéssy [3]. All (mixed) moments converge by Theorem 2.6.

We have, for simplicity, considered only the case $m/n \to \lambda > 0$, but the results are easily extended to the cases $m/n \to 0$ and $m/n \to \infty$ (at appropriate rates). Moreover, we can study the case of different probabilities $p_1, \ldots, p_n$ for the boxes by the same method; this generalization is studied in several of the references listed above. We can, further, as in some of the references, study sums $\sum_k h(N_k, k)$ for other functions $h$, possibly also depending on the box $k$. In this way, new proofs of several results in the above references may be obtained, but we leave these extensions to the reader.

## 4. Vertex degrees in random graphs

Let $\mathcal{G}_n$ be the set of all $2^{\binom{n}{2}}$ graphs with the $n$ (labelled) vertices $1, \ldots, n$. Two basic and widely studied models of random graphs are known as $G(n, m)$ and $G(n, p)$. $G(n, m)$, where $0 \leq m \leq \binom{n}{2}$, is obtained by choosing an element of $\mathcal{G}_n$ with exactly $m$ edges at random (uniformly). $G(n, p)$, where $0 \leq p \leq 1$, is defined by making a random choice for each pair of distinct vertices and connecting them by an edge with probability $p$, independently of all other edges. Note that the number of edges in $G(n, p)$ is $\mathrm{Bi}(\binom{n}{2}, p)$ and that $G(n, m)$ can be obtained as $G(n, p)$ conditioned on having exactly $m$ edges, for any $m$ and $p \in (0, 1)$. See, further, Bollobás [5] or Janson, Łuczak and Ruciński [19].

It is well known that the two random graph models $G(n, p)$ and $G(n, m)$ are very similar and, for many properties and quantities, they show the same asymptotic behaviour (for appropriate $p = p(n)$ and $m = m(n)$). In general, it is usually easy to obtain results for $G(n,p)$ from the corresponding results for $G(n,m)$, but it is often more difficult to go in the opposite direction.

For monotone properties or quantities, the situation is simple and it is possible to draw conclusions in both directions. This is well known for thresholds of monotone functions and for convergence in probability of monotone quantities. Theorem 2.3 shows that, under very general conditions, this also holds for asymptotic normality, although the asymptotic variances will usually be different for $G(n,p)$ and $G(n,m)$.

The fact that the asymptotic variances generally differ is easily seen by proceeding in the other direction, from $G(n,m)$ to $G(n,p)$; see Pittel [24]. This is a standard analysis of variance argument; the variance for $G(n,p)$ will have an extra term that can be interpreted as the part of the variance that is explained by the variation in the number of edges. In many situations, this term is of the same order as the variance for $G(n,m)$ and then the two models will have variances that are different, but of the same order. In other cases, one term may dominate the other. If the extra term is dominated by the variance for $G(n,m)$, then $G(n,m)$ and $G(n,p)$ have the same asymptotic variance. In the opposite case, the asymptotic variance for $G(n,m)$ is of a smaller order than for $G(n,p)$. It is easily seen that the latter case is exactly the case when the "first projection method" applies for $G(n,p)$; see [19], Section 6.4. In this case, it is not possible to derive precise results for $G(n,m)$ from the limit results for $G(n,p)$ (at least not without more detailed information). A typical case where our approach to $G(n,m)$ thus fails is the number of copies of a given small subgraph $H$ in $G(n,p)$ with constant $p$ and in $G(n,m)$ with $m = p\binom{n}{2}$; if $H$ has $v$ vertices, the variance is, in general, of the order $n^{2v-2}$ for $G(n,p)$ and $n^{2v-3}$ for $G(n,m)$ [15].

As an application of Nerman's theorem, we consider the numbers of vertices of different degrees in $G(n,p)$ and $G(n,m)$, with $p \sim c/n$ and $m \sim cn/2$. For $G(n,p)$, it is known that these numbers have asymptotic normal distributions; this easily extends to their joint distribution. We then use Nerman's theorem to find the same property for $G(n,m)$, but with somewhat different asymptotic variances and covariances. Recall the notation (3.3).

**Theorem 4.1.** (i) *Consider $G(n,p)$, where $p = p(n) = \lambda_n/n$ and $\lambda_n \to \lambda > 0$, and let $N_k = N_k(n)$ be the number of vertices of degree $k$, $k \geq 0$. Then,*

$$n^{-1/2}(N_k - \pi_{\lambda_n}(k)n) \xrightarrow{d} U_k, \qquad k \geq 0, \tag{4.1}$$

*jointly for all $k$, with $U_k$ jointly normal with $\mathbb{E}\, U_k = 0$ and covariances*

$$\mathrm{Cov}(U_j, U_k) = \sigma_{jk} := \pi_\lambda(j)\pi_\lambda(k)\left(\frac{(j-\lambda)(k-\lambda)}{\lambda} - 1\right) + \pi_\lambda(k)\delta_{jk}.$$

*More generally, for any sequence $(a_n)_0^\infty$ of real numbers with $a_n = O(A^n)$ for some $A < \infty$,*

$$n^{-1/2}\left(\sum_k a_k N_k - \sum_k a_k \pi_{\lambda_n}(k) n\right) \xrightarrow{d} \sum_k a_k U_k, \tag{4.2}$$

which is a normal random variable with mean 0 and variance $\sum_{j,k} a_j a_k \sigma_{jk}$.

(ii) *The same results hold for $G(n, m)$, where $m = m(n)$ and $\lambda_n := 2m/n \to \lambda > 0$, except that $U_k$ is now replaced by $\widetilde{U}_k$ with $\mathbb{E}\widetilde{U}_k = 0$ and*

$$\mathrm{Cov}(\widetilde{U}_j, \widetilde{U}_k) = \pi_\lambda(j)\pi_\lambda(k)\left(-\frac{(j-\lambda)(k-\lambda)}{\lambda} - 1\right) + \pi_\lambda(k)\delta_{jk}. \tag{4.3}$$

**Proof.** First, consider $G(n, \lambda_n/n)$. It is shown in [2], by Stein's method (see also [19], Example 6.35), that each $N_k$ is asymptotically normal. More precisely,

$$n^{-1/2}(N_k - \mathbb{E}N_k) \xrightarrow{d} U_k \sim N(0, \sigma_{kk}), \tag{4.4}$$

where

$$\sigma_{kk} := \lim_{n\to\infty} n^{-1}\mathrm{Var}\, N_k = \pi_\lambda(k)^2 \left(\frac{(k-\lambda)^2}{\lambda} - 1\right) + \pi_\lambda(k). \tag{4.5}$$

Moreover, with $p_n := \lambda_n/n$, uniformly in $k \geq 0$, we have

$$\begin{aligned}\mathbb{E}N_k &= n\binom{n-1}{k}p_n^k(1-p_n)^{n-1-k} \\ &= n\frac{\lambda_n^k}{k!}e^{-\lambda_n}\left(1 + O\left(\frac{(k+1)^2}{n}\right)\right).\end{aligned} \tag{4.6}$$

Hence, we may replace $\mathbb{E}N_k$ by $\pi_{\lambda_n}(k)n$ in (4.4).

The proof immediately extends to finite linear combinations of $N_k$, which shows joint convergence in (4.4) for all $r \geq 0$; the covariances are given by

$$\sigma_{jk} := \mathrm{Cov}(U_j, U_k) = \pi_\lambda(j)\pi_\lambda(k)\left(\frac{(j-\lambda)(k-\lambda)}{\lambda} - 1\right) + \pi_\lambda(j)\delta_{jk}.$$

Now, assume that $a_k = O(A^k)$ are given real numbers. We claim that for any such $a_k$, we have

$$n^{-1/2}\sum_{k=0}^\infty a_k(N_k - \mathbb{E}N_k) \xrightarrow{d} \sum_{k=0}^\infty a_k U_k. \tag{4.7}$$

Indeed, by the joint convergence in (4.4), this holds for the partial sums $\sum_{k=0}^K$ for any finite $K$. Moreover, from the exact formula for $\mathrm{Var}\, N_k$ (see [19], Example 6.35), it easily follows that, for any given $B$ such that $\sup_n \lambda_n \leq B$,

$$n^{-1}\mathrm{Var}\, N_k = O(B^k/k!)$$

uniformly in $k$ and it is then routine to let $K \to \infty$ to obtain (4.7) (cf. [4], Theorem 4.2). Further, again using (4.6), it follows that

$$n^{-1/2} \sum_{k=0}^{\infty} a_k(N_k - \pi_{\lambda_n}(k)n) \xrightarrow{d} \sum_{k=0}^{\infty} a_k U_k. \qquad (4.8)$$

Note, further, that by the same argument, or by the Cramér–Wold device, this holds jointly for any finite set of sequences $(a_k)$ with $a_k = O(A^k)$. In particular, since the number of edges is $M = \frac{1}{2} \sum k N_k$ and since $\sum k \pi_{\lambda_n}(k) = \lambda_n$, we have

$$n^{-1/2}(M - \tfrac{1}{2}\lambda_n n) \xrightarrow{d} V := \tfrac{1}{2} \sum_{k=0}^{\infty} k U_k. \qquad (4.9)$$

We can now transfer this result to $G(n,m)$, with $\lambda_n = 2m/n \to \lambda > 0$, as assumed in (ii). First, for any $J \geq 0$, Corollary 2.5 applies to $X_n = (N_1, \ldots, N_J)$ and $Y_n = M$, with $T$ as in (3.2). Hence, we can use Theorem 2.3, with $y_n = d_n = m = \lambda_n n/2$, and conclude from (4.1) and (4.9), which by the above hold jointly, that (4.1) also holds (jointly in all $k \geq 0$) for $G(n,m)$, with $U_k$ replaced by $\widetilde{U}_k := U_k - (\mathrm{Cov}(U_k, V)/\mathrm{Var}(V))V$; a simple calculation yields (4.3).

Moreover, a sum $\sum_k a_k U_k$ with $a_k$ increasing is stochastically increasing with respect to the number of edges and it further follows from Theorem 2.3 that (4.8) also holds for $G(n,m)$, with $U_k$ replaced by $\widetilde{U}_k$, provided $a_k$ is increasing and $a_k = O(A^k)$, again with joint convergence for several such sequences. However, any sequence $a_k = O(A^k)$ is the difference of two increasing such sequences and thus this result extends to all sequences $a_k = O(A^k)$ by linearity, which shows that (4.2) also holds for $G(n,m)$, with $U_k$ replaced by $\widetilde{U}_k$. □

**Remark 4.2.** It follows from (4.4), (4.5) and (4.6) that the mean and variances converge in (4.1) for $G(n,p)$. Hence, by Theorem 2.6, they also converge for $G(n,m)$ (with $U_k$ replaced by $\widetilde{U}_k$).

**Remark 4.3.** The limit result in (ii) is the same as for the random allocations in Section 3 with $2m$ balls. Indeed, we can add the edges to $G(n,m)$ one by one, at random, and then the vertex degrees can be described by an allocation model, as in Section 3 with $2m$ balls, except that the balls are now thrown in pairs and we condition on no pair being thrown into the same box and no two pairs being thrown into the same two boxes. Our results show, unsurprisingly, that this conditioning does not affect the asymptotic distribution of the edges.

**Remark 4.4.** The same argument applies to many other monotone functions of random graphs, for example, the size of the largest component. In this case, Pittel [24] proved asymptotic normality for both $G(n,p)$ and $G(n,m)$, with $p = c/n$ and $m = cn/2$; he first proved the result for $G(n,m)$ and then obtained the result for $G(n,p)$ as a consequence.

Nerman's theorem shows that it is also possible to do the opposite, obtaining the result for $G(n,m)$ from a result for $G(n,p)$, perhaps obtained as suggested in [16].

Another example is the size of the $k$-core, where asymptotic normality is proved in [18] for both $G(n,p)$ and $G(n,m)$, using Theorem 4.1 above. Again, Nerman's theorem shows that it suffices to study $G(n,p)$, provided we verify joint convergence with the number of edges, which provides an alternative way of treating $G(n,m)$.

## 5. Another application: spacings

To illustrate the versatility of Nerman's theorem, we give another simple application, where the random variables $Y_n$ are continuous.

Consider the $n$ spacings $S_1, \ldots, S_n$ created by $n-1$ i.i.d. random uniformly distributed points on $(0,1)$ or by $n$ such points on a circle of circumference 1. Let $a$ be a fixed positive number and let $N_a$ be the number of spacings greater than $a/n$.

If $T_1, \ldots, T_n$ are i.i.d. $\mathrm{Exp}(1)$ random variables, then

$$(S_1, \ldots, S_n) \stackrel{d}{=} \left( (T_1/n, \ldots, T_n/n) \Big| \sum_{i=1}^n T_i = n \right)$$

and thus

$$N_a \stackrel{d}{=} \left( \sum_{i=1}^n \mathbf{1}[T_i > a] \,\Big|\, \sum_{i=1}^n T_i = n \right).$$

The central limit theorem yields

$$n^{-1/2}\left( \sum_{i=1}^n \mathbf{1}[T_i > a] - n\mathrm{e}^{-a}, \sum_{i=1}^n T_i - n \right) \stackrel{d}{\longrightarrow} (X,Y),$$

where $(X,Y)$ is normal with $\mathbb{E}X = \mathbb{E}Y = 0$ and $\sigma_X^2 = \mathrm{e}^{-a}(1-\mathrm{e}^{-a})$, $\sigma_{XY} = a\mathrm{e}^{-a}$, $\sigma_Y^2 = 1$. Corollary 2.4 yields

$$n^{-1/2}(N_a - n\mathrm{e}^{-a}) \stackrel{d}{\longrightarrow} N(0, \mathrm{e}^{-a} - \mathrm{e}^{-2a} - a^2\mathrm{e}^{-2a}).$$

Moment convergence holds by Theorem 2.6.

This is a simple case of a theorem by Le Cam [21], where more general sums of the form $\sum_i h(nS_i)$ are treated. The method above applies to all monotone functions $h$ such that $\mathbb{E}h(T_1)^2 < \infty$ and more general functions can be treated by taking linear combinations. We leave the details to the reader. We can similarly study sums of the type $\sum_i h(nS_i, \ldots, nS_{i+m-1})$ and obtain, for suitable $h$, a new proof of the asymptotic normality proved by Holst [10].

## 6. A one-sided Cramér–Wold theorem

The following theorem is a version of the Cramér–Wold theorem [6], [4], Theorem 7.7; note that the Cramér–Wold theorem assumes (6.1) for arbitrary $t_1, \ldots, t_d$, but our version only assumes this when $t_i \geq 0$, which, for example, might be important for applications to monotone functions.

Recall that a random vector is *determined by its moments* if all mixed moments are finite and every random vector with the same mixed moments has the same distribution. We will use this assumption in our theorem. Note that the Cramér–Wold theorem uses no such assumption, but Example 6.4 below shows that some such condition is necessary for our one-sided version.

**Theorem 6.1.** *Suppose that $X^{(n)} = (X_1^{(n)}, \ldots, X_d^{(n)})$, $n \geq 1$, and $X = (X_1, \ldots, X_d)$ are random vectors in $\mathbb{R}^d$, where $d \geq 1$, such that*

$$\sum_{i=1}^{d} t_i X_i^{(n)} \xrightarrow{d} \sum_{i=1}^{d} t_i X_i \tag{6.1}$$

*for all real numbers $t_1, \ldots, t_d \geq 0$. Suppose, further, that the distribution of $X$ is determined by its moments. Then, $X^{(n)} \xrightarrow{d} X$. (Hence, (6.1) holds for all real $t_1, \ldots, t_d$.)*

*Remark 6.2.* A simple sufficient condition for $X$ to be determined by its moments is that $\mathbb{E} e^{a|X_k|} < \infty$ for every $k$ and some $a > 0$.

**Proof of Theorem 6.1.** Taking $t_i = \delta_{ik}$ in (6.1), we see that $X_k^{(n)} \xrightarrow{d} X_k$ for every $k$. In particular, the sequence $(X_k^{(n)})_n$ is tight for each $k$ and thus the sequence $(X^{(n)})_n$ of random vectors is tight. Consequently, every subsequence has a subsubsequence that converges in distribution and to show $X^{(n)} \xrightarrow{d} X$, it suffices to show that if some subsequence converges in distribution to $Y$, then $Y \stackrel{d}{=} X$.

Hence, assume that $Y = (Y_1, \ldots, Y_d)$ is such that $X^{(n)} \xrightarrow{d} Y$ along some subsequence. If $t_1, \ldots, t_d \geq 0$, then $\sum_k t_k X_k^{(n)} \xrightarrow{d} \sum_k t_k Y_k$ along the same subsequence and (6.1) shows that

$$\sum_{k=1}^{d} t_k Y_k \stackrel{d}{=} \sum_{k=1}^{d} t_k X_k, \qquad t_1, \ldots, t_d \geq 0. \tag{6.2}$$

In particular, with $t = (t_1, \ldots, t_d)$, denoting the characteristic function of a random vector $Z$ by $\varphi_Z(t) := \mathbb{E} e^{it \cdot Z}$, we have

$$\varphi_Y(t_1, \ldots, t_d) = \mathbb{E} e^{it \cdot Y} = \mathbb{E} e^{it \cdot X} = \varphi_X(t_1, \ldots, t_d), \qquad t_1, \ldots, t_d \geq 0.$$

The result thus follows from the following lemma. □

**Lemma 6.3.** *Suppose that $X = (X_1, \ldots, X_d)$ and $Y = (Y_1, \ldots, Y_d)$ are random vectors in $\mathbb{R}^d$, where $d \geq 1$, such that*

$$\varphi_X(t_1, \ldots, t_d) = \varphi_Y(t_1, \ldots, t_d), \qquad t_1, \ldots, t_d \geq 0. \tag{6.3}$$

*Suppose, further, that the distribution of $X$ is determined by its moments. Then, $X \stackrel{d}{=} Y$. (Equivalently, (6.3) holds for all real $t_1, \ldots, t_d$.)*

An equivalent statement is that if (6.2) holds and $X$ is determined by its moments, then $X \stackrel{d}{=} Y$ and thus (6.2) holds for all real $t_1, \ldots, t_d$. The first octant may here be replaced by any other cone with non-empty interior, using a linear change of variables.

**Proof of Lemma 6.3.** Let $1 \leq k \leq n$. Taking $t_i = 0$ for $i \neq k$, we see that $\varphi_{X_k}(t) = \varphi_{Y_k}(t)$ for $t \geq 0$ and thus for all real $t$, because $\varphi(-t) = \overline{\varphi(t)}$. Thus, $Y_k \stackrel{d}{=} X_k$ and $\mathbb{E}|Y_k|^m = \mathbb{E}|X_k|^m < \infty$ for every $m \geq 0$. It follows that $\mathbb{E}|X|^m < \infty$ and $\mathbb{E}|Y|^m < \infty$ for every $m \geq 0$. Hence, both $\varphi_X$ and $\varphi_Y$ are infinitely differentiable in $\mathbb{R}^d$. For any multi-index $\alpha$, (6.3) implies that $D^\alpha \varphi_X(t_1, \ldots, t_d) = D^\alpha \varphi_Y(t_1, \ldots, t_d)$ when $t_1, \ldots, t_d > 0$ and thus, by continuity, also when $t_1 = \cdots = t_d = 0$. Consequently,

$$\mathbb{E} X^\alpha = \mathrm{i}^{-|\alpha|} D^\alpha \varphi_X(0, \ldots, 0) = \mathrm{i}^{-|\alpha|} D^\alpha \varphi_Y(0, \ldots, 0) = \mathbb{E} Y^\alpha.$$

Thus, $X$ and $Y$ have the same moments and hence $X \stackrel{d}{=} Y$. □

Lemma 6.3 says that if a random vector is determined by its moments, its characteristic function is determined by its restriction to the first octant. The following example shows that this does not hold for all random vectors (even in two dimensions) without the extra condition that $X$ be determined by its moments. This extra condition in Theorem 6.1 and Lemma 6.3 can presumably be weakened and we leave it as an open problem to investigate more fully when a characteristic function is determined by its restriction to the first octant. (The above proof shows that Theorem 6.1 holds for such random vectors $X$.)

*Example 6.4.* Let $U, V$ and $W$ be independent random variables such that their characteristic functions $\varphi_U, \varphi_V, \varphi_W$ satisfy $\varphi_U(t) = 0$ for $|t| > 1$ and $\varphi_V(t) = \varphi_W(t)$ for $|t| \leq 1$, but $V$ and $W$ do not have the same distribution; see, for example, Feller [8], Sections XV.2 and XV.2a, for examples of such random variables. Define $X = (U+V, U-V)$ and $Y = (U+W, U-W)$. Then,

$$\varphi_X(t_1, t_2) = \varphi_U(t_1 + t_2) \varphi_V(t_1 - t_2),$$
$$\varphi_Y(t_1, t_2) = \varphi_U(t_1 + t_2) \varphi_W(t_1 - t_2).$$

If $t_1, t_2 \geq 0$, then either $t_1 + t_2 > 1$ and $\varphi_U(t_1 + t_2) = 0$, or $|t_1 - t_2| \leq t_1 + t_2 \leq 1$ and $\varphi_V(t_1 - t_2) = \varphi_W(t_1 - t_2)$; in both cases, $\varphi_X(t_1, t_2) = \varphi_Y(t_1, t_2)$, so (6.3) holds. Nevertheless, $X$ and $Y$ do not have the same distribution since $X_1 - X_2 = 2V$ and $Y_1 - Y_2 = 2W$ have different distributions. Hence, Lemma 6.3 is not true without the extra condition.

It further follows that if $t_1, t_2 \geq 0$, then $t_1 X_1 + t_2 X_2$ and $t_1 Y_1 + t_2 Y_2$ have the same characteristic function and thus the same distribution. Define $X^{(n)} = Y$ for all finite $n$. It follows that (6.1) holds (with equality for all $n$) for $t_1, t_2 \geq 0$, but $X^{(n)} \overset{d}{\not\to} X$. (Alternatively, let $X^{(n)} = Y$ for odd $n$ and $X^{(n)} = X$ for even $n$.) This shows that Theorem 6.1 also fails without the extra condition.

## Acknowledgements

I wish to thank Allan Gut, Lars Holst and Bengt Rosén for helpful comments.